\documentclass[aps, 12pt, eqsecnum]{revtex4}
\usepackage{graphicx}
\usepackage{hyperref}
\usepackage{footnote}
\usepackage{epstopdf}
\usepackage{setspace}
\usepackage{amsmath}
\usepackage{amsthm}
\usepackage{latexsym, amssymb}
\usepackage{verbatim}
\usepackage{lscape}
\usepackage{rotating}
\usepackage{float}
\usepackage[export]{adjustbox}
\usepackage{subfloat}
\usepackage{xcolor}

\def\beq{\begin{equation}}
\def\eeq#1{\label{#1}\end{equation}}
\def\eeqn{\end{equation}}

\newenvironment{Eqnarray}%
   {\arraycolsep 0.14em\begin{eqnarray}}{\end{eqnarray}}
\def\beqa{\begin{Eqnarray}}
\def\eeqa#1{\label{#1}\end{Eqnarray}}
\def\eeqan{\end{Eqnarray}}
\def\CR{\nonumber \\ }

\def\leqn#1{(\ref{#1})}

\def\ben{\begin{enumerate}}
\def\een{\end{enumerate}}

\def\lsim{\mathrel{\raise.3ex\hbox{$<$\kern-.75em\lower1ex\hbox{$\sim$}}}}
\def\gsim{\mathrel{\raise.3ex\hbox{$>$\kern-.75em\lower1ex\hbox{$\sim$}}}}

\def\Dslash{\not{\hbox{\kern-4pt $D$}}}
\def\dslash{\not{\hbox{\kern-2pt $\del$}}}
\def\over#1:#2{\frac{#1}{#2}}

\makeatletter
\def\section{\@startsection{section}{0}{\z@}{5.5ex plus .5ex minus
 1.5ex}{2.3ex plus .2ex}{\large\bf}}
\def\subsection{\@startsection{subsection}{1}{\z@}{3.5ex plus .5ex minus
 1.5ex}{1.3ex plus .2ex}{\normalsize\bf}}
\def\subsubsection{\@startsection{subsubsection}{2}{\z@}{-3.5ex plus
-1ex minus  -.2ex}{2.3ex plus .2ex}{\normalsize\sl}}

\renewcommand{\@makecaption}[2]{%
   \vskip 10pt
   \setbox\@tempboxa\hbox{\small #1: #2}
   \ifdim \wd\@tempboxa >\hsize     
       \small #1: #2\par          
     \else                        
       \hbox to\hsize{\hfil\box\@tempboxa\hfil}
   \fi}

 \def\citenum#1{{\def\@cite##1##2{##1}\cite{#1}}}
 
\newcount\@tempcntc
\def\@citex[#1]#2{\if@filesw\immediate\write\@auxout{\string\citation{#2}}\fi
  \@tempcnta\z@\@tempcntb\m@ne\def\@citea{}\@cite{\@for\@citeb:=#2\do
    {\@ifundefined
       {b@\@citeb}{\@citeo\@tempcntb\m@ne\@citea\def\@citea{,}{\bf ?}\@warning
       {Citation `\@citeb' on page \thepage \space undefined}}%
    {\setbox\z@\hbox{\global\@tempcntc0\csname b@\@citeb\endcsname\relax}%
     \ifnum\@tempcntc=\z@ \@citeo\@tempcntb\m@ne
       \@citea\def\@citea{,}\hbox{\csname b@\@citeb\endcsname}%
     \else
      \advance\@tempcntb\@ne
      \ifnum\@tempcntb=\@tempcntc
      \else\advance\@tempcntb\m@ne\@citeo
      \@tempcnta\@tempcntc\@tempcntb\@tempcntc\fi\fi}}\@citeo}{#1}}
\def\@citeo{\ifnum\@tempcnta>\@tempcntb\else\@citea\def\@citea{,}%
  \ifnum\@tempcnta=\@tempcntb\the\@tempcnta\else
  {\advance\@tempcnta\@ne\ifnum\@tempcnta=\@tempcntb \else\def\@citea{--}\fi
    \advance\@tempcnta\m@ne\the\@tempcnta\@citea\the\@tempcntb}\fi\fi}
\makeatother

\newenvironment{Abstract}{\begin{quotation} \begin{center}
                       ABSTRACT
     \end{center}\bigskip  }{\end{quotation}}

\def\Acknowledgements{\bigskip  \bigskip \begin{center} \begin{large}
             \bf ACKNOWLEDGEMENTS \end{large}\end{center}}

\newcommand{\dm}{\mathrm{dim\,} }

\newcommand{\bs}{\backslash}

\theoremstyle{plain}
\newtheorem{thm}{Theorem}[section]
\newtheorem{lem}{Lemma}[section]

\newtheorem{cor}{Corollary}[section]

\theoremstyle{definition}
\newtheorem{defn}{Definition}[section]

\newtheorem{exmp}{Example}[section]

\theoremstyle{remark}
\newtheorem*{rem}{Remark}

\begin{document}

\begin{titlepage}
\title{The Inductive Graph Dimension from The Minimum Edge Clique Cover}
\author{Kassahun Betre$^1$, Evatt Salinger$^2$}
\affiliation{1 San Jose State University, 1 Washington Square, Sci 148, San Jose, CA 95192\\
2 Pepperdine University, 24255 Pacific Coast Hwy, Malibu CA 90263, USA}

\maketitle

\begin{Abstract}
In this paper we prove that the inductively defined graph dimension has a simple additive property under the join operation. The dimension of the join of two simple graphs is one plus the sum of the dimensions of the component graphs: $\dm (G_1+ G_2) = 1 +\dm G_1+ \dm G_2$. We use this formula to derive an expression for the inductive dimension of an arbitrary finite simple graph from its minimum edge clique cover. A corollary of the formula is that any arbitrary finite simple graph whose maximal cliques are all of order $N$ has dimension $N-1$. We finish by finding lower and upper bounds on the inductive dimension of a simple graph in terms of its clique number. 

\end{Abstract}
\vfill
\vfill
\end{titlepage}
\def\thefootnote{\fnsymbol{footnote}}
\setcounter{footnote}{0}
\tableofcontents
\newpage
\setcounter{page}{1}

\section{Introduction}
In classical and quantum physics, the dimension of space is a given starting quantity. But in quantum gravity that is not as clear cut. Spacetime is not a concept that holds together at sub-Planckian lengths. In theories of dynamically generated geometric space from pre-geometric combinatorial spaces, such as background independent formulations of quantum gravity, it becomes paramount to find the pre-geometric combinatorial definition of dimension that maps to the dimension of the emergent space. In general, one may want to define the dimension of discrete combinatorial spaces, such as finite graphs or abstract simplicial complexes. One may look to dimension theory to gain insight into this, however, the effort will come out empty handed in the field of dimension theory of topological spaces. 

Dimension theory of topological spaces had developed into a mature field by the late 1960's \cite{pears}. The rigorous exposition of dimension theory was started by Brouwer in 1913 and independently by Menger and Urysohn in 1922  \cite{brouwer, HurewiczWallman}. They gave a definition of dimension for a separable metric topological space that matches with the intuitive dimension for Euclidean space. The Menger-Urysohn dimension is defined inductively and is a topological invariant. It assigns dimensions of $-1$ to the empty set and a positive semi-definite integer to any complete metric space. In particular, any finite or countably infinite space has dimension 0.  

In an attempt to assign dimension to finite simple graphs, one may restrict oneself to connected graphs and think of them as metric spaces with the natural metric that assigns a length between two vertices equal to the number of edges in the shortest path between them. This metric turns the graph into a topological space with the discrete topology (through the usual definition of the open ball). If thought of as a topological space, one can then apply the Menger-Urysohn inductive dimension on the graphs. However this gives the trivial dimension of zero to all finite simple graphs by virtue of the simple fact that they are finite spaces. 
Another way to assign dimension to finite simple graphs is to first construct an abstract simplicial complex $K$ from the graph, then consider the geometric realization $|K|$ of the complex and assign its dimension as a topological space to be the dimension of the graph. Depending on the choice of $K$, one again gets either a somewhat trivial result of always zero or one (if $K$ consists only of 0-simplices for vertices and 1-simplices for edges), or some integer value (if $K$ is the Whitney complex). The later assignment, while more informative, still washes a lot of important connectivity features of the graph and lumps many disparate graphs into the same dimension. That is not to mention the unnecessary intermediate step involving geometric realizations. 

Other notions of dimension such as fractal dimensions also lead to trivial results of zero to all finite graphs since the fractal dimensions are given as limits. An attempt to generalize the Hausdorff fractal dimension to finite simple graphs is presented by Alonso in \cite{2016arXiv160708130A}. However, this graph dimension definition suffers from a highly unnatural assignment of infinity to all complete graphs $K_N$ in contrast to the intuitive and natural notion that the dimension of the tetrahedron ($K_4$) should be $3$. 

To our knowledge, a purely combinatorial (graph-theoretic) and natural definition of the dimension of a finite simple graph is one given by Oliver Knill in \cite{2010arXiv1009.2292K} and explored in greater depth in \cite{2011arXiv1112.5749K}. Motivated by the Menger-Urysohn inductive dimension of topological spaces, this graph dimension is defined inductively. It exploits the inherent connectivity structure of graphs as it is defined based on induced unit spheres, and satisfies some important topological properties expected from dimensions. For example, it is shown in \cite{Knill:2015ku} that the Cartesian product of simple graphs obeys the Kuenneth formula $\dm(G_1\times G_2) \le \dm G_1 + \dm G_2$. Though the inductive graph dimension is not a topological invariant under the discrete topology, it admits a subbasis for a topology whereby homeomorphic graphs have the same set of dimensions of their subbasis \cite{knill2014notion}.

An important construction in algebraic topology with direct parallel in graph theory is the notion of the {\em join} also known in graph theory as the Zykov sum \cite{zykov, hararay, kozlov}. In algebraic topology, the Urysohn-Menger inductive dimension of the geometric realization of two complexes is 1 plus the sum of the dimensions of each. For example, the join of two one-simplices is the tetrahedron which is a 3-simplex. One natural question to ask of the Knill inductive graph dimension and the join operation on graphs is whether or not a similar relationship holds, i.e., whether or not the dimension of the join of two graphs is one plus the sum of the dimension of each. This is in fact shown to be the case for the restricted class of geometric graphs (basically vertex schemes of combinatorial manifolds). 
In \cite{Knill:2015ku}, the inductive graph dimension is shown to obey $$\dm \left(S^{k}+ S^{l}\right) = 1+ k + l,$$ where $S^{k}$ is a $k-$dimensional geometric sphere, defined recursively as a graph such that every unit sphere in the graph is a $(k-1)-$ dimensional geometric sphere with the base case $S^0$ being a disconnected graph with two vertices. In this paper we will show that this formula generalizes to the join of any finite simple graphs, i.e., that for any two finite simple graphs $G_1$ and $G_2$, the inductive dimension of the join of the two graphs is 
\beq 
\dm (G_1+ G_2) = 1 + \dm G_1 + \dm G_2.
\eeq{eq:dimSum} 
This general formula will be used to derive an expression for the inductive dimension of a graph from its minimum clique cover, and to place bounds on the inductive dimension. 



We begin section 2 with a review of the graph theory definitions relevant to our discussion. We introduce a few non-standard notations for convenience, such as $K_V$ for the complete graph formed by connecting each vertex in the set $V$ by an edge, and the minimum and maximum clique number of a graph. Section 3 gives the proof of \leqn{eq:dimSum}. In section 4 we derive a formula for the inductive dimension of a graph from the minimum clique cover. In section 5 we place some bounds on the inductive dimension.

\section{The Inductive Dimension}
\begin{defn}


The complete graph over $N$ vertices is $K_N$. Let us define an analogous concept. If $V$ is a set of vertices, we define $K_V$ to be the complete graph whose vertex set is $V$. For example, $K_{\{1,2,3\}}$ is the complete graph with vertex set $V = \{1,2,3\}$.
\end{defn}
	
	The order of the biggest maximal clique of a graph $G$ is called the {\bf clique number} of the graph and denoted by $\omega(G)$. This is a standard definition and notation of the clique number of a graph in the literature. However, we will find it convenient to separate the clique number in two; the {\bf maximum clique number} referring to the order of the largest maximal clique (denoted by $\omega(G)$), and the {\bf minimum clique number} referring to the order of the smallest maximal clique. There is no standard notation in the literature for the minimum clique number, so we will use the notation $\gamma(G)$ to denote the order of the smallest maximal clique in $G$.

\begin{defn}

The {\bf edge clique cover} ECC of a graph is a set of cliques that cover the edges of the graph; i.e, the union of the edge sets of the cliques in the edge clique cover is equal to the edge set of the graph. The minimum edge clique cover uses the fewest possible number of cliques, so the cliques in the minimum edge clique cover are maximal cliques. The minimum number of cliques needed to form an edge clique cover is the edge clique cover number, denoted by $\theta_e(G)$. However, the edge clique cover need not cover all vertices of the graph; for example, isolated vertices of the graph are not contained by any clique in the edge clique cover. 

In this paper we will only be concerned with the edge clique cover number (ECC) and not the vertex clique cover (VCC) and therefore we will drop the explicit mention of ``edge" and simply refer to the clique cover. 
\end{defn}

\begin{defn}
A graph is {\bf pure} if all maximal cliques have the same order.
\end{defn}

\begin{defn}
	A {\bf sphere} in $G$ of radius $r$ centered at the vertex $v$, denoted by $S_G(v, r)$ is an induced subgraph of $G$ whose vertices are all vertices in $G$ with distance of exactly $r$ from $v$. 

The unit sphere will have $r=1$. We will abbreviate the notation for the unit sphere at $v$ and write simply $S_G(v)$ without the $1$. The radius is assumed to be one unless explicitly stated. The vertices of the unit sphere at $v$ are the {\bf neighbors} of $v$. The degree of a vertex is therefore equal to the order of the unit sphere at the vertex.	
\end{defn}

\begin{defn}
	A {\bf ball} in $G$ of radius $r$ centered at $v$, denoted by $B_G(v,r)$, is the induced subgraph whose vertices are all vertices in $G$ at a distance of $r$ or less from $v$. We write simply $B_G(v)$ for the unit ball with $r=1$.
\end{defn}


\begin{rem}
	The unit ball $B_G(v)$ is the join of the vertex $v$ with the unit sphere $S_G(v)$. 
			\beq
				B_G(v) = v + S_G(v)
			\eeqn
\end{rem}

\begin{defn}
	The {\bf inductive dimension} of a graph, denoted by $\text{dim}(G)$, is defined recursively as 
	\beqa
		\dm(G) &=& \begin{cases}
		-1,&\quad \text{if } G \text{ is the empty graph}\\
		\frac{1}{|G|}\sum_v \dm_G(v), &\quad \text{otherwise}\\
		\end{cases}\CR\CR
		\dm_G(v) &=& 1 + \dm S_G(v)
	\eeqa{eq:dimension}
\end{defn}

\begin{rem}
		The completely disconnected graph over $N$ nodes has inductive dimension 0.
		The complete graph over $N$ nodes $K_N$ has dimension of $N-1$. The inductive dimension matches the Euclidean dimension of simplexes if the complete graph $K_N$ is thought of as a $(N-1)-$simplex embedded in $N$-dimensional Euclidean space.
		All bipartite graphs are 1-dimensional. That is because $S_G(v)$ is completely isolated graph for any $v$ in a bipartite graph; so that $\dm_G(v) = 1 + \dm S_G(v) = 1$, and $\dm G = \frac{1}{|G|}\sum_v 1 = 1.$ All tree graphs are 1-dimensional since all tree graphs are bipartite.  All Cycle graphs $C_n$ are 1-dimensional for $n>3$. 
\end{rem}
	
\section{Dimension of The Join of Graphs}

One immediate observation of the inductive dimension given in \cite{2011arXiv1112.5749K} is that the dimension of a disconnected graph with two components is the weighted average of the dimensions of the components. This is true for any disconnected graphs with arbitrary number of components.  We will restate this observation in its generality in the lemma below.

\begin{lem}
\label{lem:DisconnectedDim}
Let $G$ be disjoint union of components $G_1, G_2, \dots, G_N$, then 
	\beq
		\dm G =  \sum_{i=1}^N \frac{|G_i|}{|G|} \dm G_i.
	\eeqn
\end{lem}

{\it Proof:}
We begin by observing that for a vertex $v \in G_i, \;S_G(v) = S_{G_i}(v)$, therefore $\dm_G(v) = \dm_{G_i}(v)$. Then, 
\beqa
|G|\dm G &=& \sum_{v\in G} \dm_G(v)\CR
&=&\sum_{i=1}^N\sum_{v\in G_i} \dm_{G}(v)\CR
&=&\sum_{i=1}^N\sum_{v\in G_i} \dm_{G_i}(v)\CR
&=&\sum_{i=1}^N|G_i|\dm G_i\nonumber
\eeqan

\begin{lem}
\label{lem:dimsum}
The dimension of the join of two arbitrary graphs is one plus the sum of the dimensions of each; i.e.,  
\beq
	\dm (G_1 + G_2) = 1+ \dm G_1 + \dm G_2. 
\eeq{eq:dimsim}

\end{lem}

\noindent{\it Proof:}
This lemma generalizes a corollary 7 of \cite{2011arXiv1112.5749K}, which states the join of a $k-$dimensional and an $\ell-$dimensional geometric spheres is a geometric sphere of dimension $n = k + \ell + 1$. We will prove the theorem by induction over the order of $G_1 + G_2$.
The base case is when either $|G_1| = 0$ or $|G_2| = 0$. Since the dimension of the empty graph is by definition $-1$, we have 
\beqa
\dm (G_1 + \emptyset) = 1 + \dm G_1 + \dm \emptyset = \dm G_1.\nonumber
\eeqan

Assume the formula is true for arbitrary $|G_1+ G_2| = |G_1| + |G_2| = k$ and neither graph is empty, we will show that the formula holds when the order of the join is $k + 1$. Let's call the join $G_1 +  G_2 = G$ with $|G| = |G_1| + |G_2| = k+1$. First we observe that if a vertex $v \in G_1$, then $S_G(v) = S_{G_1}(v) +  G_2$ and similarly if $v\in G_2$, then $S_G(v) = S_{G_2}(v) +  G_1$ (See Fig. \ref{fig:G1PlusG2} below). Then,
\beqa
|G|\dm(G) &=& \sum_{v\in G_1}\dm_{G}(v) + \sum_{v\in G_2}\dm_{G}(v)\CR
	(k+1)\dm(G) &=& \sum_{v\in G_1}\left(1 + \dm S_{G}(v)\right) + \sum_{v\in G_2} \left(1 + \dm S_{G}(v)\right)\CR
	&=& \sum_{v\in G_1}\Big(1 + \dm \left(S_{G_1}(v)+ G_2\right)\Big) + \sum_{v\in G_2} \Big(1 + \dm \left(G_1+ S_{G_2}(v)\right)\Big)\CR
	&=& \sum_{v\in G_1}\Big(1 + (1 + \dm S_{G_1}(v)+ \dm G_2)\Big) + \sum_{v\in G_2} \Big(1 + (1 + \dm G_1+\dm S_{G_2}(v))\Big)\CR
	&=& |G_1| + |G_1|\dm G_1 + |G_1|\dm G_2 + |G_2| + |G_2|\dm G_1 + |G_2|\dm G_2\CR
	&=& \left( |G_1| + |G_2| \right)\left(1 + \dm G_1 + \dm G_2 \right)\CR
	\dm (G_1 +  G_2) &=& 1 + \dm G_1 + \dm G_2\nonumber
\eeqan
In the fourth line we applied the inductive assumption that the formula holds when the order of the join graph is $k$, which applies here since $ |S_{G_1}(v) +  G_2| \leq k,$ and similarly $|G_1+ S_{G_2}(v)| \leq k$ for any $v$. 
\begin{figure}
\vspace{-0.5cm}
\begin{center}
  \includegraphics[width=0.8\linewidth]{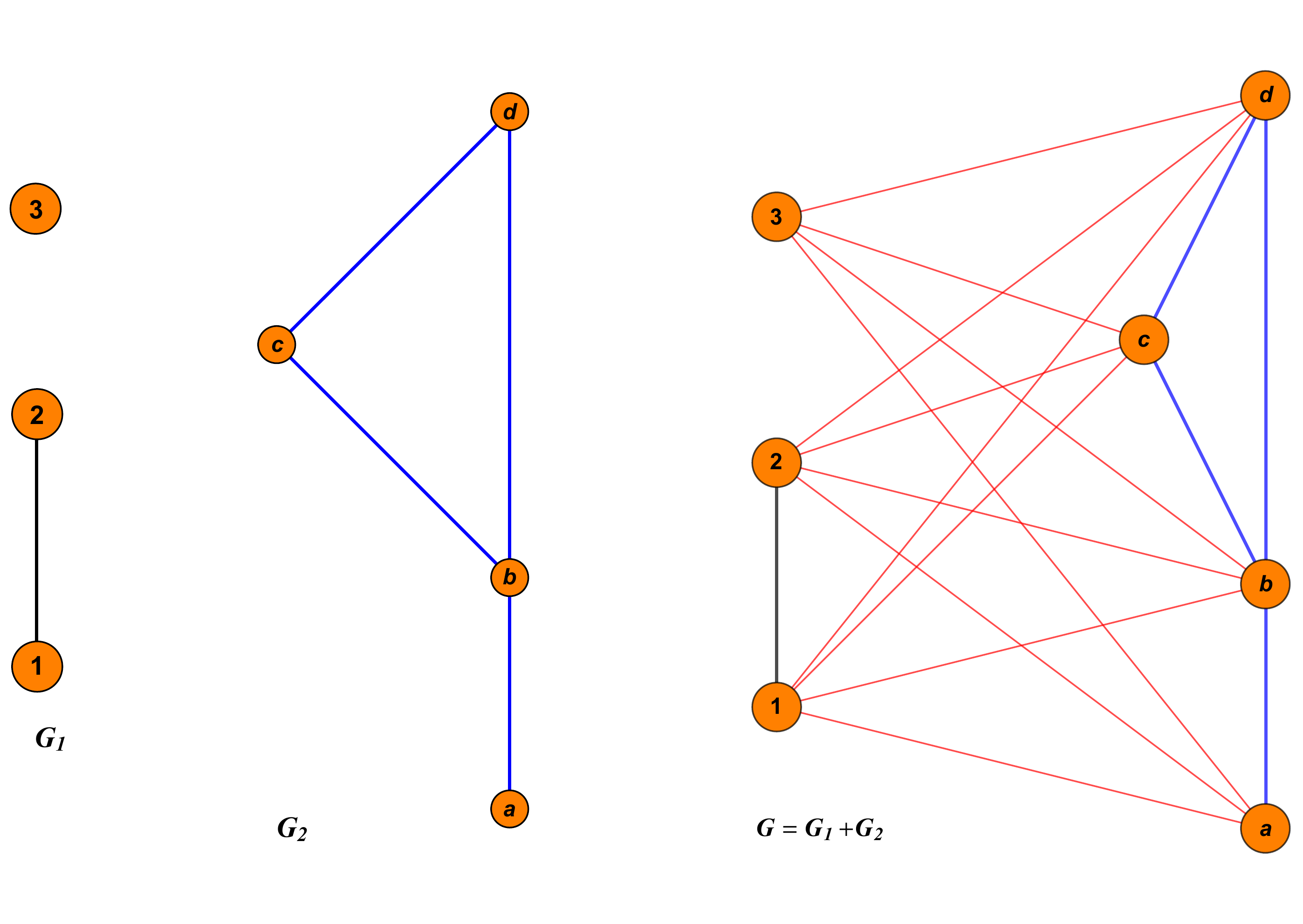}
  \caption{Two graphs $G_1$ and $G_2$ are shown on the left. Their join is on the right. For any vertex in $G_1$, the sphere in $G$ at $v$ is the join of the sphere in $G_1$ at $v$ and $G_2$. The same holds for any vertex $v \in G_2 \implies S_G(v) = S_{G_2}(v) +  G_1.$}
 \label{fig:G1PlusG2}
\end{center}
\end{figure}
\begin{cor}
If $G = G_1 +  G_2 +  \dots +  G_k$, then $\dm G = (k-1)+\dm G_1 + \dm G_2 + \dots +\dm G_k.$
\end{cor}
{\it Proof:}
\beqa
\dm G &=& \dm \Big(G_1 +   (G_2 +  \dots +  G_k)\Big)\CR
&=& 1 + \dm G_1 + \dm \left(G_2 +  G_3 +  \dots +  G_k \right) \CR
&=& 2 + \dm G_1 + \dm G_2 + \dm \left(G_3 +  G_4 +  \dots +  G_k \right)\CR
\vdots\CR
&=& (k-1) + \dm G_1 + \dm G_2  + \dots + \dm G_k
\eeqan

\begin{exmp}
The complete graph $K_N$ is the join of $N$ copies of $K_1$. So, $\dm K_N = \dm (K_1 +  K_1 +  \dots +  K_1) = (N-1) + N \dm K_1 = N-1.$
\end{exmp}

\begin{cor}
\label{cor:DimBallEqualsOnePlusDimSphere}
The dimension of any vertex in $G$ is the dimension of the unit ball at that vertex.
	\beq
		\dm_G (v) = \dm B_G(v) = 1 + \dm S_G(v).
	\eeqn
\end{cor}

{\it Proof:}
This follows immediately from the fact that $B_G(v) = v +  S_G(v)$ so that $\dm _G(v) = 1 + \dm S_G(v) = \dm \Big(v +  S_G(v)\Big) = \dm B_G(v)$.

\section{The inductive dimension from the Minimum Edge Clique Cover}

In this section, let us assume that the minimum edge clique cover of the graph $G$ is known. Let $m$ be the clique cover number, and let $K_{V_1},K_{V_2},\dots,K_{V_m}$ be the maximal cliques in the minimum clique cover with vertex sets $V_1, V_2, \dots, V_m$. Let us introduce some notations that will simplify the discussion. We will use the notation $K_{V_iV_j}$ to denote $G[V_i\cap V_j] $, the induced graph in $G$ over the vertices in the intersection of $K_{V_i}$ and $K_{V_j}$. By extension, we will use $K_{V_iV_jV_k} = G[V_i\cap V_j\cap V_k]$, etc. Since the intersection of two complete graphs is another complete graph, each of the graphs $K_{V_iV_j}, K_{V_iV_jV_k}, \dots$ is also a complete graph. We will also use $G[K_{V_i}\cup K_{V_j}]$ to refer to the induced graph $G[V_i \cup V_j]$ over the union of the vertex sets of $K_{V_i}$ and $K_{V_j}$. However, this graph is in general not a complete graph.

It will be useful to define a notation for the number of vertices in a given maximal clique that are not contained in any of the other maximal cliques in the clique cover. Let $||K_{V_i}||$ be the number of vertices in $K_{V_i}$ that are not in any of the other cliques. Similarly, let $||K_{V_iV_j}||$ be the number of vertices contained in the intersection $K_{V_iV_j}$ but not in any of the other maximal cliques, etc. Then, by the inclusion-exclusion principle, 

	\beqa
		||K_{V_i}|| &=& |K_{V_i}| - \sum_{\substack{j\neq i\\j=1}}^{m}|K_{V_iV_j}|+\sum_{\substack{j,k\neq i\\j<k = 1}}^{m}|K_{V_iV_jV_k}|-\dots \pm |K_{V_1V_2\dots V_m}|\CR
		||K_{V_iV_j}|| &=& |K_{V_iV_j}| - \sum_{\substack{k\neq i,j\\k=1}}^{m}|K_{V_iV_jV_k}| +\sum_{\substack{k_1,k_2\neq i,j\\k_1<k_2 = 1}}^{m}|K_{V_iV_jV_{k_1}V_{k_2}}|-\dots \pm |K_{V_1V_2\dots V_m}|\CR
		\vdots\CR	
		||K_{V_1V_2\dots V_m}|| &=& |K_{V_1V_2\dots V_m}|	
	\eeqan
Therefore, 
	\beqa
		|G| &=& \sum_{i=1}^{m}||K_{V_i}|| + \sum_{i>j=1}^m||K_{V_iV_j}||+ \sum_{i_1>i_2>i_3=1}^m||K_{V_{i_1}V_{i_2}V_{i_3}}|| + \dots\CR
		&&+ \sum_{i_1>i_2>\dots>i_{m-1}=1}^m||K_{V_{i_1}V_{i_2}\dots V_{i_{m-1}}}|| \CR
		&&+ ||K_{V_1V_2\dots V_m}||\nonumber
	\eeqa{eq:OrderOfGfromCliques}
In words this is simply saying that the vertices in $G$ can be partitioned into those contained in a single maximal clique, the intersection of only two maximal cliques, etc.

\begin{lem}
\label{lem:TwoMaxCliques}
Let G be a graph with clique cover number of 2 so that $G = K_{V_1} \cup K_{V_2}$. Then, 
	\beq
		\Big(|G| - |K_{V_1V_2}|\Big)\dm G = ||K_{V_1}||\dm K_{V_1} + ||K_{V_2}||\dm K_{V_2}
	\eeqn  
\end{lem}
{\it Proof:}
If we remove the $|K_{V_1V_2}|$ vertices in $G$ we end up with a disconnected graph $G^0 = G - K_{V_1V_2}$ with two complete graph components, $K_{V_1\bs (V_1\cap V_2)}$ and $K_{V_2\bs (V_1\cap V_2)}$. The order of the resulting graph is $|G^0| = |G| - |V_{V_1V_2}| = ||K_{V_1}|| + ||K_{V_2}||.$ Using lemma \ref{lem:DisconnectedDim},
\beqa
|G^0|\dm G^0 &=& |K_{V_1\bs (V_1\cap V_2)}|\dm \left(K_{V_1\bs (V_1\cap V_2)} \right) + |K_{V_1\bs (V_1\cap V_2)}|\dm \left(K_{V_1\bs (V_1\cap V_2)} \right)\CR
&=& ||K_{V_1}||(|K_{V_1}| - |K_{V_1V_2}|) + ||K_{V_2}||(|K_{V_2}| - |K_{V_1V_2}|)\nonumber
\eeqan
\beqa
|G^0|\dm G^0 &=& ||K_{V_1}||\dm K_{V_1} + ||K_{V_2}||\dm K_{V_2} - \dm K_{V_1V_2}(||K_{V_1}|| + ||K_{V_2}||)\CR
&=&  ||K_{V_1}||\dm K_{V_1} + ||K_{V_2}||\dm K_{V_2} -|G^0| \dm K_{V_1V_2}
\eeqan
Then, note that the graph $G$ is the join of $G^0$ and $K_{V_1V_2}$, $G = G^0 + K_{V_1V_2}$, and so using lemma \ref{lem:dimsum},
\beqa
\dm G &=& 1 + \dm G^0 + \dm K_{V_1V_2}\CR
&=& \dm G^0 + |K_{V_1V_2}|\CR
|G^0|\dm G &=& |G^0|\dm G^0 + |G^0| |K_{V_1V_2}| \CR
&=&  ||K_{V_1}||\dm K_{V_1} + ||K_{V_2}||\dm K_{V_2} \CR
\Big(|G| - |K_{V_1V_2}| \Big)\dm G &=&||K_{V_1}||\dm K_{V_1} + ||K_{V_2}||\dm K_{V_2}\nonumber
\eeqan
\begin{thm}
\label{thm:KnillDimFromCliques}
Let $\{K_{V_1}, K_{V_2}, \dots,K_{V_m}\}$ be the set of maximal cliques in the minimum clique cover of $G$ and let $K_L = \bigcap_{i=1}^m K_{V_i}$. Then, 
\beqa
	\Big(|G| - |K_L|\Big)\dm G =&&  \sum_{i=1}^{m} ||K_{V_i}||\; \dm K_{V_i}\\
	&+&\sum_{i>j=1}^{m} ||K_{V_iV_j}||\; \dm G[K_{V_i} \cup K_{V_j}]\CR
		&+&\sum_{i_1>i_2>i_3=1}^{m} ||K_{V_iV_jV_k}||\; \dm G[K_{V_i} \cup K_{V_j}\cup K_{V_k}]\CR
		&+&\vdots\CR
		&+&\sum_{i_1>i_2>\dots>i_{m-1}}^{m}||K_{V_{i_1}V_{i_2}\dots V_{i_{m-1}}}||\; \dm G[\bigcup_{i\in\{i_1,i_2,\dots,i_{m-1}\}}K_{V_i}].\nonumber
\eeqa{eq:KnillDimCliqueFormula} 
\end{thm}
\noindent The sums $\sum_{i_1>i_2>\dots>i_k=1}^{m}$ etc. are sums over all $k$ subsets of the set $\{1,2,\dots,m\}$. 
The theorem is stating that to find the dimension of a graph with $m$ cliques in the minimum clique cover, we can proceed vertex by vertex, first selecting all vertices that are contained in only one clique, whose dimension in $G$ is simply the dimension of the cliques that contain them, then those vertices in the intersection of only two cliques, etc.  

{\it Proof:}
We will do inductive proof on the clique number $m$. The case when $m=1$ is trivial. The case $m=2$ is proven in lemma \ref{lem:TwoMaxCliques}. 

For the inductive step, assume the formula in theorem \ref{thm:KnillDimFromCliques} holds true for $m$ maximal cliques in the minimum clique cover of $G$. For the case when $G$ has $m+1$ maximal cliques, let $K_L = \bigcap_{i=1}^{m+1}K_{V_i}$, and let $G^0 = G - K_L$. Then, if $K_{V_i}$ is a maximal clique in $G$, then $K_{V_i\bs V(K_L)} = K_{V_i} - K_L$ is a maximal clique in $G^0$. So, the cliques $K_{V_1}-K_L, K_{V_2}-K_L, \dots, K_{V_{m+1}}-K_L$ form the maximal cliques in a minimum clique cover of $G^0$.

Further, note that the numbers $||K_{V_i}||, ||K_{V_iV_j}||, \dots$ are the same in $G^0$ and $G$ since the vertices removed in $G^0$ are contained in the intersection of all maximal cliques and so are counted only in $K_L$. Then,
\beqa
|G^0|\dm G^0 &=& \sum_{i=1}^{m+1} ||K_{V_i}||\; \dm G[K_{V_i} - K_L]\\
		+&&\sum_{i>j=1}^{m+1} ||K_{V_iV_j}||\; \dm G[\left(K_{V_i}- K_L\right) \cup \left(K_{V_j}- K_L\right)]\CR
		+&&\sum_{i>j>k=1}^{m+1} ||K_{V_iV_jV_k}||\; \dm G[\left(K_{V_i}- K_L\right) \cup \left( K_{V_j}- K_L\right)\cup \left( K_{V_j}- K_L\right)]\CR
		+&&\vdots\CR
		+&&\sum_{i_1>i_2>\dots>i_{m}=1}^{m+1}||K_{V_{i_1}V_{i_2}\dots V_{i_{m}}}||\; \dm G[\bigcup_{i\in\{i_1,i_2,\dots,i_{m}\}}\left(K_{V_i}- K_L\right)]\nonumber
\eeqan
For each of the induced subgraphs, we have 
\beqa
	\dm G[\left( K_{V_{i_1}} - K_L\right) \cup \dots \cup\left( K_{V_{i_k}}-K_L\right)] &=& \dm G[\left( K_{V_{i_1}} \cup \dots \cup K_{V_{i_k}}\right)-K_L]\CR
	&=& \dm G[K_{V_{i_1}} \cup \dots \cup K_{V_{i_k}}]  - |K_L|.\nonumber
\eeqan
That is because the induced graph $G[K_{V_{i_1}} \cup \dots \cup K_{V_{i_k}}]$ is the join of $G[\left( K_{V_{i_1}} - K_L\right) \cup \dots \cup\left( K_{V_{i_k}}-K_L\right)]$ and $K_L$, and by the inductive hypothesis (since the induced graphs all have clique number of $m$ or less). Therefore,
\beqa
|G^0|\dm G^0 =&& \Bigg(\sum_{i=1}^{m+1} ||K_{V_i}||\; \dm G[K_{V_i}]\\
		&& +\sum_{i>j=1}^{m+1} ||K_{V_iV_j}||\; \dm G[K_{V_i} \cup K_{V_j}]\CR
		&& +\vdots\CR
		&& +\sum_{i_1>i_2>\dots>i_{m}=1}^{m+1} ||K_{V_{i_1}V_{i_2}\dots V_{i_{m}}}||\; \dm G[\bigcup_{i\in\{i_1,i_2,\dots,i_{m}\}}K_{V_i}]\Bigg) - |K_L||G^0|\nonumber
\eeqan
In the last term of the last line we have used \leqn{eq:OrderOfGfromCliques}.

Coming back to $G$, since $G = G^0 + K_L$, we again have $\dm G = \dm G^0 + |K_L|$, and
\beqa
|G^0|\dm G &=& |G^0|\dm G^0 + |G^0||K_L|\CR
&=&  \sum_{i=1}^{m+1} ||K_{V_i}||\; \dm G[K_{V_i}]\\
		&& +\sum_{i>j=1}^{m+1} ||K_{V_iV_j}||\; \dm G[K_{V_i} \cup K_{V_j}]\CR
		&& +\sum_{i>j>k=1}^{m+1} ||K_{V_iV_jV_k}||\; \dm G[K_{V_i} \cup K_{V_j}\cup K_{V_j}]\CR
		&& +\vdots\CR
		&& +\sum_{i_1>i_2\dots>i_{m}=1}^{m+1} ||K_{V_{i_1}V_{i_2}\dots V_{i_{m-1}}}||\; \dm G[\bigcup_{i\in\{i_1,i_2,\dots,i_{m}\}}K_{V_i}]\nonumber
\eeqan
This completes the proof of theorem \ref{thm:KnillDimFromCliques}.

This theorem allows us to show in the corollary below that the dimension of a pure graph is equal to the dimension of the maximal cliques.  

\begin{cor}
\label{cor:RegularCliques}
Let $K_{V_1}, K_{V_2}, \dots, K_{V_m}$ be the maximal cliques in the minimum clique cover of a pure graph $G$. If $|V_1| = |V_2| = \dots, |V_m| = N$ then,  $\dm G = \dm K_N = N-1$. Furthermore, every vertex in the graph has equal dimension of $N-1$.
\end{cor}

{\it Proof:}
We will do an inductive proof on the clique number $m$. Let $K_L = \bigcap_{i=1}^{m}K_{V_i}$. The case of $m=2$ follows immediately from the formula for $\dm G$ in lemma \ref{lem:TwoMaxCliques}.
\beqa
	\Big(|G|-|K_L| \Big)\dm G &=& ||K_{V_1}||\dm K_N + ||K_{V_2}||\dm K_N\CR
	&=&\dm K_N\Big( ||K_{V_1}|| + ||K_{V_2}||\Big) = \dm K_N\Big( |G| - |K_L|\Big)\CR
	\dm G &=& \dm K_N = N-1\nonumber
\eeqan

Assume the case is true for $m$ maximal cliques, all of the same order $N$, so that $\dm G = \dm K_N$. Then, since the dimension of the graph with $m+1$ maximal cliques is expressed as the sum of dimensions of graphs each with a clique cover number of $m$ or less, each dimension in the sum will be the same ($\dm K_N$) by the inductive hypothesis. It then follows  
\beqa
\Big(|G|-|K_L| \Big)\dm G &=& \dm K_N\Bigg(\sum_{i=1}^{m}||K_{V_i}|| + \sum_{i>j=1}^{m}||K_{V_iV_j}||+\dots\CR
&&\qquad\qquad\quad+\sum_{i_1>i_2>\dots>i_{m-1}=1}^{m}||K_{V_{i_1}V_{i_2}\dots V_{i_{m-1}}}||\Bigg)\CR
&=&\dim K_N\Big(|G| - |K_L| \Big)\CR
\dm G &=& \dm K_N.\nonumber
\eeqan

To prove that every vertex in the graph has regular dimension $N-1$, we note that the sphere at any vertex in the graph is the union of complete graphs of uniform order $N-1$. Therefore the sphere at each vertex has dimension $N-2$, giving the result that the dimension at each vertex is $N-1$.

\begin{figure}[H]
\vspace{-0.5cm}
\begin{center}
  \includegraphics[width=0.9\linewidth]{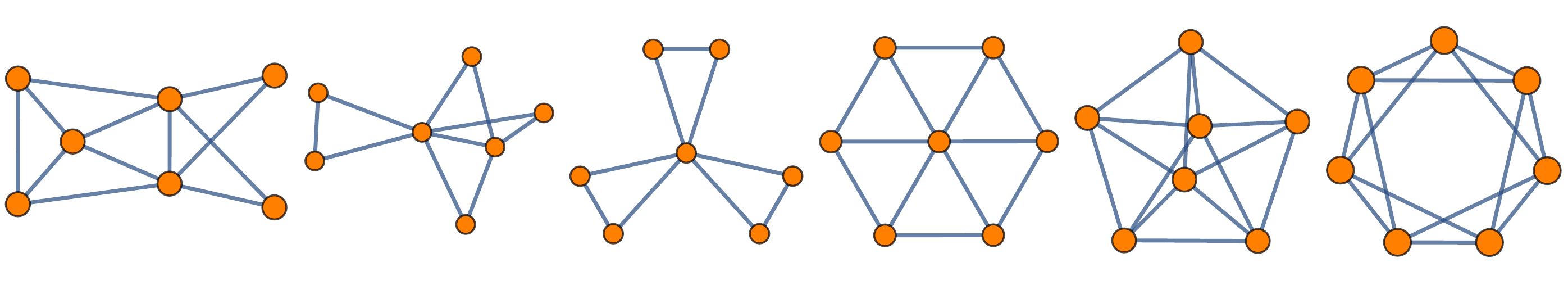}
  \caption{Examples of pure graphs with 7 vertices and $K_3$ maximal cliques. All have inductive dimension of 2.}
 \label{fig:N7Pure2DGraphs}
\end{center}
\end{figure}

\begin{rem}
The corollary gives another way to see that all trees and forests are 1-dimensional since trees and forests are pure graphs with $K_2$ maximal cliques. Furthermore, the complete $k-$partite graph has dimension $k-1$. 
\end{rem}
\begin{exmp}
Each of the graphs in Fig (\ref{fig:N7Pure2DGraphs}) are pure with $K_3$ maximal cliques. Therefore, all of them have inductive dimension of 2, and each vertex of each graph is has dimension 2. 
\end{exmp}

%
%


\section{Graphs of Uniform Dimension}
\begin{figure}[H]
\vspace{-0.0cm}
\begin{center}
  \includegraphics[width=0.5\linewidth]{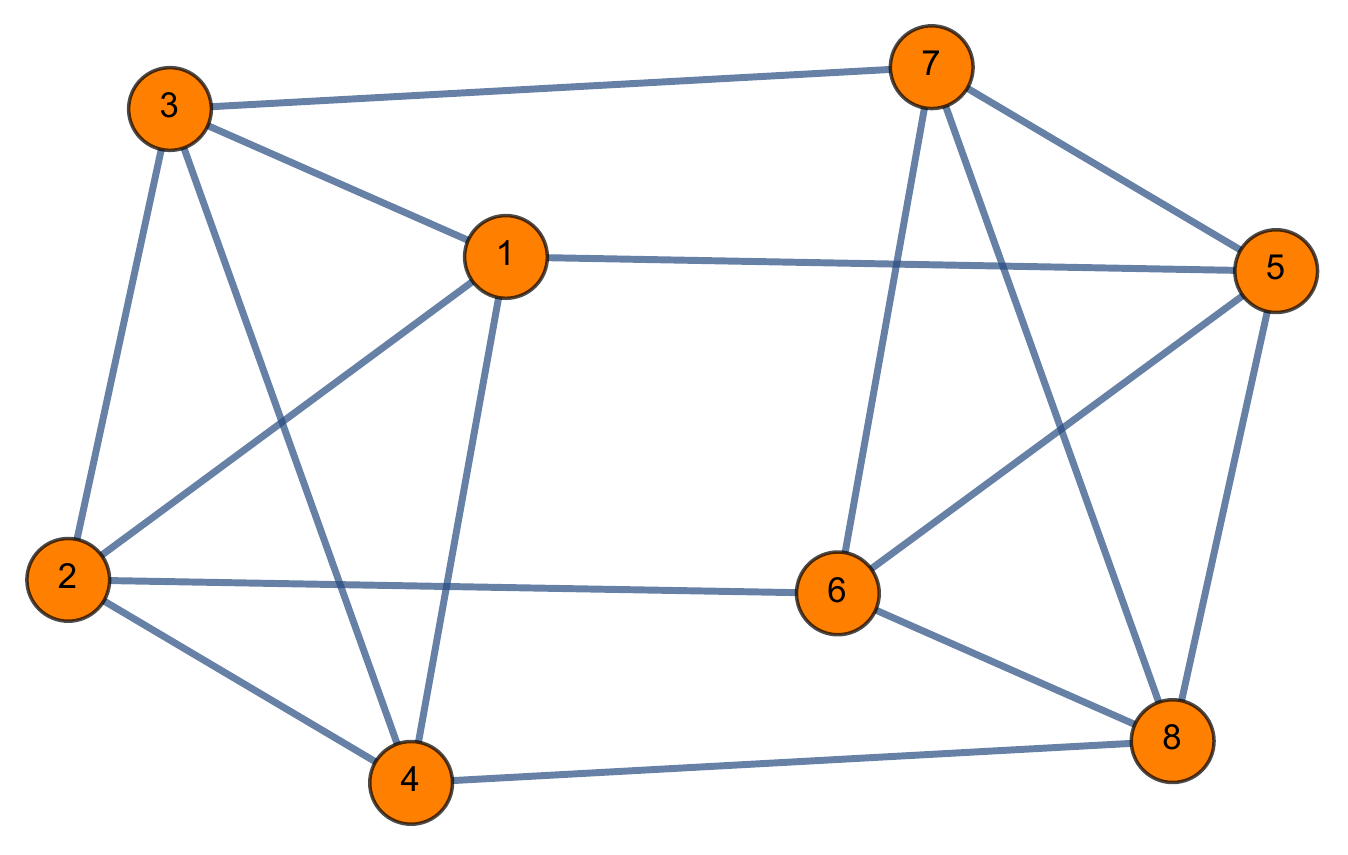}
  \caption{A regular dimensional graph that is not pure. Each vertex has dimension $5/2$. The maximal cliques are $K_4$ and $K_2$.}
 \label{fig:ImpureD5over2graph}
\end{center}
\end{figure}

\begin{figure}[H]
\vspace{-2.0cm}
\begin{center}
  \includegraphics[width=1.0\linewidth]{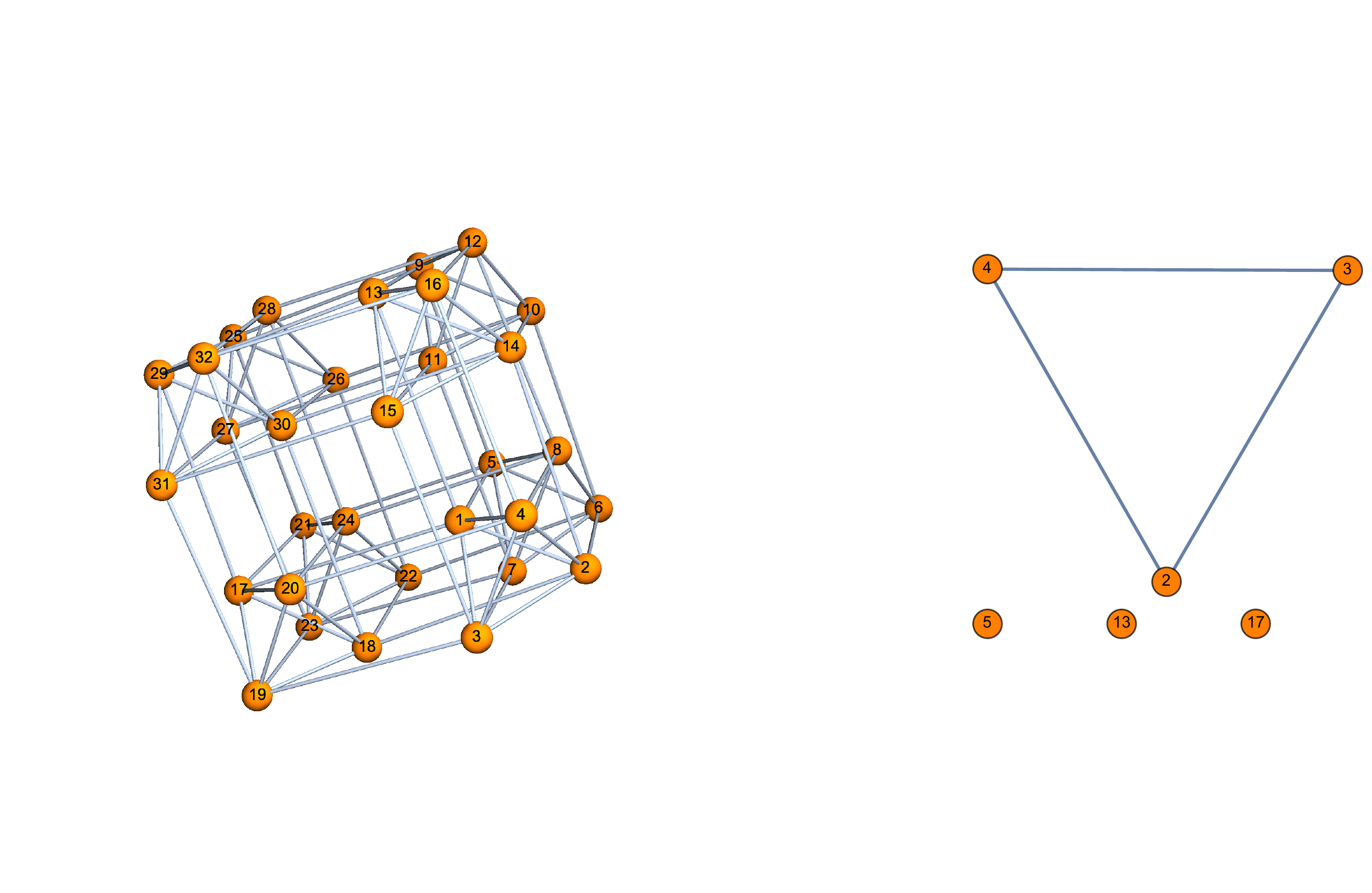}
  \vspace{-2.0cm}
  \caption{ A regular 2-dimensional graph that is not pure (Right) and the unit sphere at vertex 1 (Left). Every unit sphere at each vertex is isomorphic to graph on the right.}
 \label{fig:Impure2Dgraph}
\end{center}
\end{figure}

Let us call a graph {\bf uniform} if the dimension of each vertex is equal. While a pure graph with all maximal cliques of order $N$ has dimension $N-1$, and hence is uniform, the converse is not true. If a graph is uniform with dimension $d$ at each vertex, it is not necessarily true that the graph is pure. One reason for this is that the dimension of each vertex of a uniform graph may be a non-integer rational number. For example, the graph shown in Fig (\ref{fig:ImpureD5over2graph}) has dimension $5/2$ at each vertex. However, the graph is not pure since there are maximal cliques of order 4 and 2 intersecting at each vertex. The situation does not improve even if we restrict the dimension at each vertex to be integral. There are uniform graphs with  integral dimension $d$ at each vertex that are not pure. For example, the graph on the right side in Fig (\ref{fig:Impure2Dgraph}) has regular dimension of 2 at each vertex. However the graph is not pure since maximal cliques have orders 4 and 2. The graph is constructed by replacing each vertex of the cube by a $K_4$ and by adding a parallel edge from each vertex of the $K_4$ to a vertex in the $K_4$ that replaced the (what used to be adjacent) node. The graph is {\em homogenous}; the sphere at each vertex looks like the disconnected graph on the right side.

\section{Bounds on the Inductive Dimension}

In this section we derive some bounds on the inductive dimension. The bounds are given in terms of the maximum clique number $\omega(G)$ and the minimum clique number $\gamma(G)$ of the graph defined in section 2. Let $\gamma(G) = \ell, \omega(G) = k$. From theorem \ref{thm:KnillDimFromCliques} and corollary \ref{cor:RegularCliques} it follows that $\ell-1 \le \dm(G) \le k-1$. The lower bound comes from turning all maximal cliques in $G$ into $K_{\ell}$ by removing edges without changing the order of the graph. This results in a graph with dimension $\ell-1$. The upper bound comes from turning all maximal cliques in $G$ into complete graphs $K_k$ by adding edges, which gives $k -1$ for the the dimension of the graph.

\begin{cor}
\label{cor:dimbound1}
If $G$ is a graph with maximum clique number $\omega(G) = k$, then 
\beq
	\frac{k(k-1)}{|G|}\le \dm G \le k-1.
\eeq{eq:dimbound1}
\end{cor}

{\it Proof:}
The lower bound occurs for the edge-minimal graph with maximum clique number $k$. This is a graph with a single $k-$clique and every other vertex isolated, i.e., $G$ is the disjoint union of $K_k$ and $|G|-k$ isolated vertices. Using \ref{lem:DisconnectedDim}, $\dm G = k/|G|\,\dm K_k = k(k-1)/|G|.$ The upper bound occurs for edge maximal graphs with clique number $k$, which is when all maximal cliques in $G$ have order $k = \omega(G)$. By corollary \ref{cor:RegularCliques} the dimension of $G$ is $k-1$.

For connected graphs the lower bound in \leqn{eq:dimbound1} can be improved.

\begin{cor}
\label{cor:dimbound2}
Let $G$ be a connected graph with maximum clique number $\omega(G) = k$, then 
\beq
	1 + \frac{k^2(k-1)(k-2)}{|G|\Big( k(k-2) + |G|\Big)}\le \dm G \le k-1.
\eeq{eq:dimbound2}
\end{cor}
{\it Proof:}
Let the edge set of the maximal clique of order $k$ be $U$, so that $K_U$ is the highest order clique in the clique cover. The lower bound comes from a graph with a single maximal clique of order $k$ and all other maximal cliques of order 2; i.e., the union of a tree with $K_U$. For such a graph, every vertex that is not in $K_U$ will be one dimensional. The vertices $v\in K_U$ will have unit spheres that are the disjoint union of a complete graph of order $k-1$ and $({\rm deg}(v) - (k-1))$ isolated vertices. Therefore these unit spheres will have dimension $(k-2)(k-1)/{\rm deg}(v).$ It follows,
\beqa
|G|\dm G &=& \sum_{v\in K_U}\dm_G(v) + \sum_{v\notin K_U}\dm_G(v) = \sum_{v\in K_U}\Big(1 + \dm S_G(v) \Big)+ \sum_{v\notin K_U}1 \CR
&=&k +\left( \sum_{v\in K_U}\frac{(k-2)(k-1)}{{\rm deg}(v)}\right) + |G| - k\CR
\dm G&=& 1 + \frac{(k-2)(k-1)}{|G|}\sum_{v\in K_U}\frac{1}{{\rm deg}(v)}
\eeqan

\begin{figure}[H]
\begin{center}
\vspace{-0.5cm}
  \includegraphics[width=0.5\linewidth]{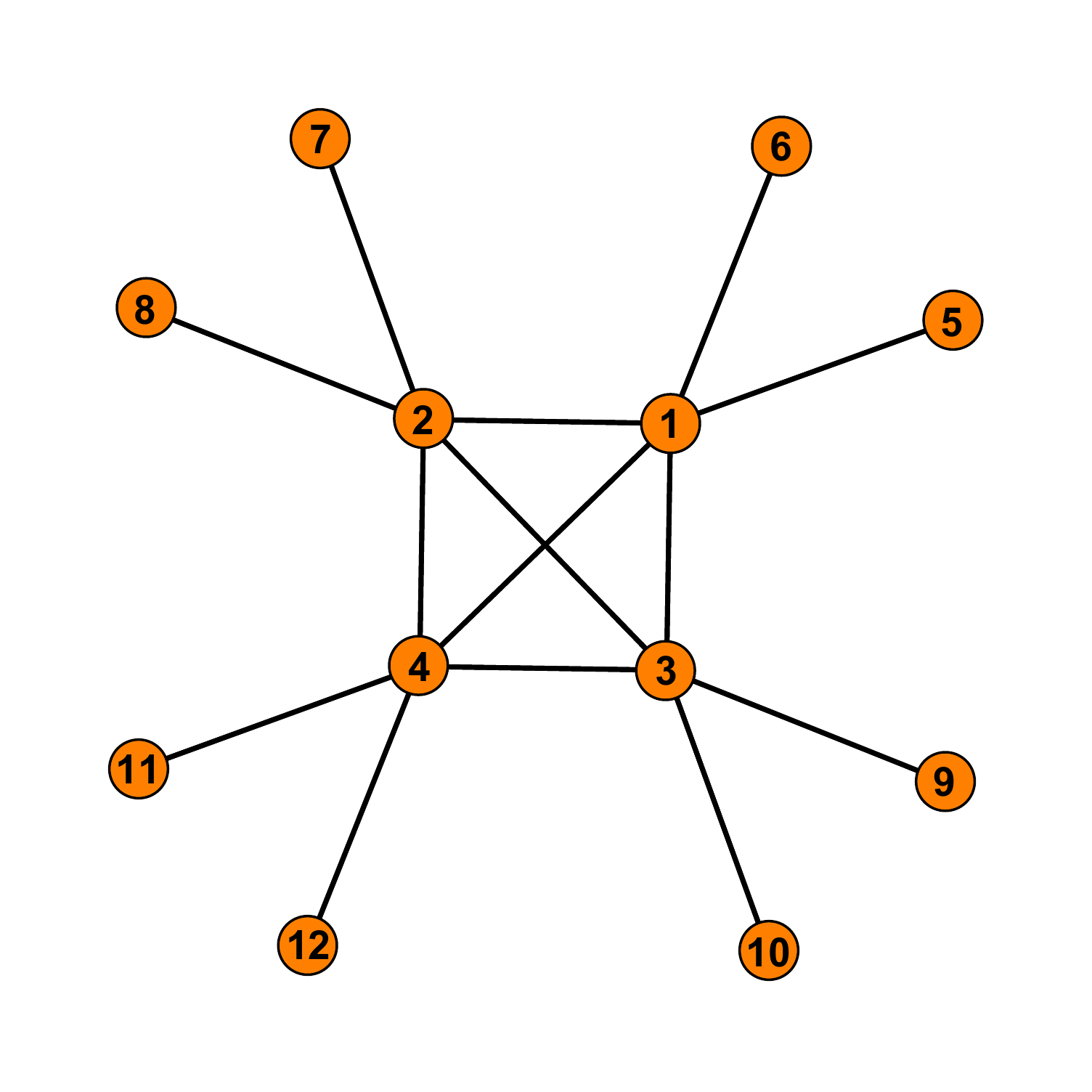}
\vspace{-0.3cm}
  \caption{A graph that saturates the lower bound of the dimension formula in \leqn{eq:dimbound2} for $|G|=12$ and $k=4$.}
 \label{fig:StarGraph}
\end{center}
\end{figure}

The sum $\sum_{v\in K_U} 1/{\rm deg}(v)$ can be minimized subject to the constraint $\sum_{v\in K_U}{\rm deg}(v) \le k\Big((k-1) + (|G|-k)/k\Big)$. The upper bound of the constraint comes from having every vertex not in $K_U$ be connected to a vertex in $K_U$ to maximize the degree of the vertices in $K_U$. Since $1/x$ is a convex function, the inverse of the mean is always less than or equal to the mean of the inverses, and the minimum value for $\sum_{v\in K_U} 1/{\rm deg}(v)$ is achieved when every term in the sum is equal to the inverse of the average degree. In that case the degree of every vertex in $K_U$ is $(k-1)+(|G|-k)/k$. This results in a lower bound of 
\beqa
\dm G&=& 1 + \frac{k^2(k-1)(k-2)}{|G|\Big( k(k-2) + |G|\Big)}.\nonumber
\eeqan
The graphs that saturate the lower bound are graphs that look like the star graph, but with the central vertex replaced by a complete graph $K_k$, and where the leaves are distributed as equitably as possible among the vertices in the complete graph. See Fig. \ref{fig:StarGraph}

\section{Summary}
\label{sec:summary}
The main results of this paper are the following:
\begin{enumerate}
	\item The dimension of the join of two graphs is the sum of their dimensions plus one.
	\beq
	\dm \left(G_1+G_2\right) = \dm G_1 + \dm G_2 + 1.\nonumber
	\eeqn
	\item The inductive dimension of a graph has the formula given in theorem \ref{thm:KnillDimFromCliques} in terms of the maximal cliques in the minimum clique cover of the graph.
	\item Any finite pure graph with all maximal cliques of the same order $N$ has dimension $N-1$. In addition,  such a graph has regular vertex dimensions of $N-1$.
	\item The dimension of a graph is bounded between $k(k-1)/|G|$ and $k-1$ where $k=\omega(G)$ is the maximum clique number. If the graph is connected, the lower bound can be improved to \leqn{eq:dimbound2},
\beq
	1 + \frac{k^2(k-1)(k-2)}{|G|\Big( k(k-2) + |G|\Big)}\le \dm G \le k-1.\nonumber
\eeqn
\end{enumerate}
\Acknowledgements
KB would like to thank Kevin Iga at Pepperdine university for simplifying the proof of lemma \ref{lem:dimsum} and for offering valuable feedback on the earlier drafts of the paper. He would also like to thank Patrick Wells at UC Davis and Oliver Knill at Harvard for the valuable conversations. ES is supported in part by AYURI fund from Pepperdine University.

%


\end{document}